\documentclass[12pt,twoside]{article}
\usepackage{amsmath,amssymb,amsthm}
\usepackage{graphicx}

\newtheorem{theorem}{Theorem}[section]
\newtheorem{proposition}[theorem]{Proposition}

\newtheorem{corollary}[theorem]{Corollary}
\newtheorem{example}[theorem]{Example}
\newtheorem{remark}[theorem]{Remark}

\newcommand{\Hom}{\operatorname{Hom}}
\newcommand{\Epm}{\operatorname{Epm}}

\title{\textbf{Jordan Left $\alpha$-centralizers on Algebras\\ with Applications to Group Algebras}}

\author{
\textbf{Mojdeh Eisaei}\\
Department of Mathematics\\
Payame Noor University\\
Shiraz, Iran\\
\texttt{mojdehessaei59@student.pnu.ac.ir}
\and
\textbf{Mohammad Javad Mehdipour}\thanks{Corresponding author}\\
Department of Mathematics\\
Shiraz University of Technology\\
Shiraz 71555-313, Iran\\
\texttt{mehdipour@sutech.ac.ir}
\and
\textbf{Gholam Reza Moghimi}\\
Department of Mathematics\\
Payame Noor University\\
Shiraz, Iran\\
\texttt{moghimimath@pnu.ac.ir}
}

\begin{document}

\maketitle

\begin{abstract}
We prove that every Jordan left $\alpha$-centralizer from an algebra $A$ with a right identity into an arbitrary algebra $B$ is a left $\alpha$-centralizer. This implies all Jordan homomorphisms between such algebras are homomorphisms. We extend this result to continuous Jordan left $\alpha$-centralizers when $A$ has a bounded left approximate identity. For the group algebra $L^1(G)$, we characterize weakly compact Jordan left $\alpha$-centralizers when $\alpha$ is continuous and surjective, showing $L^1(G)$ admits a weakly compact epimorphism if and only if $G$ is finite. Consequently, the existence of a non-zero $\alpha$-derivation on $L^1(G)$ is equivalent to $G$ being compact and non-abelian.
\end{abstract}

\section{Introduction}

Let $G$ be a locally compact group with its group algebra $L^1(G)$ and the measure algebra $M(G)$. 
It is well-known that $M(G)$ is the dual of $C_0(G)$, the space of all complex-valued continuous functions on $G$ that vanish at infinity. For any two measures $\mu, \nu \in M(G)$, their convolution product is defined by 
\[
\mu\ast \nu (f) = \int_G \int_G f(xy) d\mu(x) d \nu(y)
\]
for all $f \in C_0(G)$. With this product $M(G)$ becomes a unital Banach algebra.
We note that $L^1(G)$ is a closed ideal of $M(G)$ and consequently forms a Banach algebra.

Let $A$ and $B$ be algebras. We denote by $\Hom(A, B)$ the set of all homomorphisms from $A$ into $B$, and by $\Epm(A, B)$ the set of all epimorphisms from $A$ into $B$. When $A$ and $B$ are normed algebras, we write $\Hom_c(A, B)$ and $\Epm_c(A, B)$ for the corresponding sets of continuous homomorphisms and continuous epimorphisms, respectively. In the case where $A=B$, we simplify the notation to $\Hom(A)$, $\Hom_c(A)$, $\Epm(A)$ and $\Epm_c(A)$.

Given a homomorphism $\alpha\in\Hom(A, B)$, an additive function $T: A\rightarrow B$ is called a \textit{left $\alpha$-centralizer} if it satisfies
\[ T(ax)= T(a) \alpha(x) \]
for all $a,x \in A$. The function $T$ is called a \textit{Jordan left $\alpha$-centralizer} if it satisfies
\[ T(a^2)= T(a) \alpha(a) \]
for all $a\in A$. When $A=B$ and $\alpha$ is the identity map, $T$ is termed a \emph{left centralizer} or \emph{Jordan left centralizer}, respectively. For $\alpha\in\Hom(A, B)$ and $b_0\in B$, one can construct the canonical left $\alpha$-centralizer $\Lambda_{b_0, \alpha}: B\rightarrow B$ defined by
\[ \Lambda_{b_0, \alpha}=\lambda_{b_0}\circ\alpha, \]
where $\lambda_{b_0}: B\rightarrow B$ denotes the left $\alpha$-centralizer $\lambda_{b_0}(b)=b_0 b$.

The centralizers on various algebraic structures have been studied by several authors. For example, Wendel [17] investigated the left centralizers on group algebras and gave a characterization for $L^1(G)$. This work was continued by Akemann [2] and Sakai [16], who examined the weakly compact case and proved the fundamental equivalence between the existence of non-zero weakly compact left centralizers on $L^1(G)$ and the compactness of $G$. Subsequently, Zalar [18] studied Jordan left centralizers in the more general setting of semiprime rings and proved that every Jordan left centralizer on a semiprime ring with characteristic different from 2 is a left centralizer; see [5, 6, 7] for study of centralizers on other Banach algebras. Motivated by these foundational contributions, the present work investigates Jordan left $\alpha$-centralizers on both general algebras and group algebras, where $\alpha\in\hbox{Hom}(A, B)$.

In Section 2, we prove that for an algebra $A$ with a right identity, every Jordan left $\alpha$-centralizer from $A$ into an arbitrary algebra $B$ is a left $\alpha$-centralizer. Furthermore, we examine this result to the case where A is a Banach algebra with a bounded left approximate identity and $\alpha$ is continuous, and prove the same result for continuous Jordan left $\alpha$-centralizers. As applications of these results, we show that both Jordan homomorphisms and their continuous versions between such algebras are homomorphisms. 

Section 3 focuses on the group algebra $L^1(G)$, where $\alpha\in\hbox{Hom}_c(L^1(G))$. We characterize continuous Jordan left $\alpha$-centralizers on $L^1(G)$ and provide a necessary and sufficient condition for the existence of a non-zero weakly compact Jordan left $\alpha$-centralizer. Furthermore, we establish that the existence of a weakly compact epimorphism on $L^1(G)$ is equivalent to the finiteness of $G$. As a consequence of our results, we prove that $L^1(G)$ admits a non-zero weakly compact $\alpha$-derivation if and only if $G$ is compact and non-abelian provided that $\alpha$ is a continuous epimorphism. Finally, we give some examples to show that the surjectivity condition on $\alpha$ is essential in our results and can not be omitted.
%-----------------------------------------------------------------------
\section{\normalsize\bf Jordan left $\alpha$-centralizers on algebras} 

The main result of this section is the following.
%-----------------------------------------------------------------------
%-----------------------------------------------------------------------
%-----------------------------------------------------------------------
\begin{theorem} \label{m1} Let $A$ be an algebra with a right identity $u$, $B$ be any algebra, and $\alpha\in\emph{Hom}(A, B)$. Then every Jordan left $\alpha$-centralizer $T: A\rightarrow B$ is a 
left $\alpha$-centralizer. 
\end{theorem}
%-----------------------------------------------------------------------
%-----------------------------------------------------------------------
%-----------------------------------------------------------------------
{\it Proof.}
Let $T:A \rightarrow B$ be a Jordan left $\alpha$-centralizer. 
Then 
$T(a^2)=T(a) \alpha (a)$ 
for all $a \in A$. 
The linearity of this equation yields
\begin{eqnarray} \label{e1} 
T(ax+xa)= T(a) \alpha(x) + T(x) \alpha(a)
\end{eqnarray}
for all $a, x \in A$. This implies that 
\begin{eqnarray}\label{zei}
T(a^2x+xa^2)&=& T(a)\alpha(ax)+T(x)\alpha(a^2)\nonumber\\
&=&T(a)\alpha(ax)+T(ax+xa)\alpha(a)-T(a)\alpha(xa)
\end{eqnarray}
for all $a, x\in A$. Substituting $x=ax+xa$ into (1) and applying (2), we obtain
\begin{eqnarray} \label{e6} 
T(axa)= T(a) \alpha(x) \alpha(a). 
\end{eqnarray}
Using the right identity $u$, we observe
\begin{eqnarray} \label{e7}
T(ua)= T(u au ) &=&T(u) \alpha(a) \alpha(u) \nonumber \\
&=& T(u) \alpha (au) \\
&=& T(u) \alpha(a). \nonumber
\end{eqnarray}
For $z\in A$, define $r:=z-uz$. Then $Ar=\{0\}$. 
Taking $x= u+r$ in (3), we have
\begin{eqnarray} \label{e9}
T(ux+rx) &=& T((u+r)x(u+r))\nonumber\\
&=&T(u+r)\alpha(x)\alpha(u+r) \nonumber \\
&=& T(u+r) \alpha(x) \\
&=& T(u) \alpha(x) + T(r) \alpha(x)\nonumber
\end{eqnarray}
for all $x \in A$. On the other hand, by (4) 
\begin{eqnarray*} 
T(ux+rx) =T(ux) + T(rx)= T(u) \alpha(x)+ T(rx). 
\end{eqnarray*} 
From this and (5) we see that 
\begin{eqnarray*} 
T(rx)=T(r) \alpha(x). 
\end{eqnarray*} 
Hence for every $a, x \in A$, we have
\begin{eqnarray*} 
T(ax) &=& T(ax) - T(u) \alpha (ax) + T(u) \alpha (ax) \\
&=& T(ax) - T(uax) + T(u) \alpha(ax) \\
&=& T((a-ua)x) + T(u) \alpha(ax) \\
&=& T(a-ua) \alpha(x) + T(u) \alpha(ax) \\
&=& T(a) \alpha(x) - T(ua) \alpha(x)+ T(u) \alpha(ax) \\
&=& T(a) \alpha (x). 
\end{eqnarray*} 
Therefore, $T$ is a left $\alpha$-centralizer. 
$\hfill\square$\\

For algebras $A$ and $B$, let us recall that a linear map $\alpha: A\rightarrow B$ is called a \emph{Jordan homomorphism} if 
$$
\alpha(a^2)=\alpha(a)^2
$$
for all $a\in A$. It is clear that every Jordan homomorphism $\alpha: A\rightarrow B$ is a Jordan left $\alpha$-centralizer. So Theorem 2.1 yields the following result.

\begin{corollary} Let $A$ be an algebra with a right identity and $B$ be any algebra. Then every Jordan homomorphism from $A$ into $B$ is a homomorphism.
\end{corollary}
%-----------------------------------------------------------------------
%-----------------------------------------------------------------------
%-----------------------------------------------------------------------
\begin{remark} \label{m10}
{\rm (i) Following a proof analogous to Theorem 2.1, one can show that the same result holds for rings with a certain characteristic.

(ii) Let $A$ be a unital algebra with identity $1_A$, let $B$ be any algebra, and $\alpha\in\hbox{Hom}(A, B)$. Then every Jordan left $\alpha$-centralizer $T: A \rightarrow B$ is a left $\alpha$-centralizer and we have 
\begin{eqnarray*} 
T(a)= T(1_Aa) = T(1_A) \alpha(a)=\Lambda_{T(1_A), \alpha}(a)
\end{eqnarray*}
for all $a \in A$. So $T$ coincides with the canonical $\alpha$-centralizer $\Lambda_{a_0, \alpha}$, where $a_0=T(1_A)$. It would be valuable to examine whether this result extends to algebras endowed merely with a right identity.

(iii) Several important classes of group algebras possess either right identities or identities. For instance,
the Banach algebras $L_0^\infty(G)^*$ as defined in [10] and $L^1(G)^{**}$ admit right identities. Also, the Fourier-Stieltjes algebra $B(G)$, the measure algebra $M(G)$, the Banach algebra $(M(G)_0^*)^*$ as defined in [13], and the Banach algebras $\hbox{LUC}(G)^*$ all possess identities, where $\hbox{LUC}(G)$ is the space of all left uniformly continuous functions on $G$; for more study on these Banach algebras see [1, 4, 8, 13, 14]. Hence every Jordan left $\alpha$-centralizer on these Banach algebras is automatically a left $\alpha$-centralizer. So every Jordan homomorphism between such Banach algebras is necessarily a homomorphism.}
\end{remark}
%-----------------------------------------------------------------------
%-----------------------------------------------------------------------
%-----------------------------------------------------------------------
\begin{corollary}
Let $A$ be a unital normed algebra, $B$ be any normed algebra, and $\alpha\in\emph{Hom}_c(A, B)$. Then every Jordan left $\alpha$-centralizer $T: A\rightarrow B$ is continuous. 
\end{corollary}
%-----------------------------------------------------------------------
%-----------------------------------------------------------------------
%-----------------------------------------------------------------------
{\it Proof.}
Let $1_A$ be the identity element of $A$, and let $T: A\rightarrow B$ be a Jordan left $\alpha$-centralizer. By Remark 2.3 (ii), there exists $b_0\in B$ such that 
$$T(a)=\Lambda_{b_0, \alpha}(a)$$
for all $a \in A$. Since $\alpha$ is continuous, we have $\|\alpha(a)\|\leq\|\alpha\|\|a\|$ for all $a\in A$. Consequently, we obtain the following norm estimate. 
\begin{eqnarray*}
\| T(a) \| = \| \Lambda_{b_0, \alpha}(a) \| \leq \|b_0\| \| \alpha \| \| a\|. 
\end{eqnarray*}
This inequality shows that $T$ is continuous. 
$\hfill\square$\\

%-----------------------------------------------------------------------
%-----------------------------------------------------------------------
%-----------------------------------------------------------------------
%-----------------------------------------------------------------------
%-----------------------------------------------------------------------
%-----------------------------------------------------------------------
%-----------------------------------------------------------------------
%-----------------------------------------------------------------------
%-----------------------------------------------------------------------
%-----------------------------------------------------------------------
%-----------------------------------------------------------------------
%-----------------------------------------------------------------------
Let $A $ be a Banach algebra.
In the next result, we equip the second conjugate of $A $ with the first Arens product. 
Let us recall that for $m,n \in A ^{**}$, the first Arens product $m$ and $n$ is defined by $\langle mn, f \rangle= \langle m, nf \rangle$, where $\langle nf, a \rangle= \langle n, f a \rangle$ and $\langle fa, x \rangle = \langle f, ax \rangle$ for all $f \in A ^*$ and $a, x \in A $. 

%-----------------------------------------------------------------------
%-----------------------------------------------------------------------
%-----------------------------------------------------------------------
\begin{theorem} \label{m7}
Let $A $ be a Banach algebra with a bounded left approximate identity, $B$ be any Banach algebra, and $\alpha\in\emph{Hom}_c(A, B)$. 
Then every continuous Jordan left $\alpha$-centralizer $T: A\rightarrow B$ is a left $\alpha$-centralizer. 
\end{theorem}
%-----------------------------------------------------------------------
%-----------------------------------------------------------------------
%-----------------------------------------------------------------------
%-----------------------------------------------------------------------
{\it Proof.} 
Let $T: A \rightarrow B$ be a continuous Jordan left $\alpha$-centralizer. 
We first observe that $T^*(f) \, a = \alpha^* \, ( f \, T(a))$ for all $f \in A^*$ and $a\in A$. So for any $n \in A^{**}$, we have
$$n \, T^* (f)= T^* ( \alpha ^{**} (n) f).$$ 
This implies that $T^{**}: A^{**}\rightarrow B^{**}$ 
is a Jordan left $\alpha^{**}$-centralizer. Hence
$$T^{**} (mn) = T^{**}(m) \alpha^{**} (n)$$ for all $m,n \in A ^{**}$. 
Since $A $ possesses a bounded left approximate identity, its second dual admits a right identity. 
In view of Theorem 2.1, $T^{**}$ is a left $\alpha^{**}$-centralizer. 
Finally, noting that 
$$T^{**}|_{A}= T\quad\hbox{and}\quad\alpha^{**}|_{A}=\alpha,$$ 
we conclude that
$T$ is a left $\alpha$-centralizer. 
$\hfill\square$\\ 

%-----------------------------------------------------------------------
%-----------------------------------------------------------------------
%-----------------------------------------------------------------------
As an immediate consequence of Theorem 2.5, we have the following result.

\begin{corollary} Every continuous Jordan homomorphism on a Banach algebra with a bounded left approximate identity is a homomorphism.
\end{corollary}

Let us remark that the important classes of Banach algebras, including the group algebra $L^1(G)$, the Fourier algebra $A(G)$, any $C^*$-algebra and every amenable Banach algebra possess bounded approximate identities. Hence every continuous Jordan left $\alpha$-centralizer (respectively, homomorphism) on these Banach algebras, is a left $\alpha$-centralizer (respectively, homomorphism), where $\alpha$ is a continuous homomorphism on them.

\begin{proposition} Let $A$ be a Banach algebra with a bounded approximate identity, $B$ be any Banach algebra and $\alpha\in\emph{Hom}(A, B)$. If $T: A \rightarrow B$ is a continuous Jordan left $\alpha$-centralizer, then the following statements hold.

\emph{(i)} If $T$ is weakly compact, then there exists $b_0\in B$ such that $T=\Lambda_{b_0, \alpha}$. 

\emph{(ii)} If $\alpha$ is weakly compact and $T$ is surjective, then $\alpha=T\circ\lambda_{a_0}$ for some $a_0\in A$.
\end{proposition}
{\it Proof.}
Let $(e_i)$ be a bounded approximate identity for $A$. If $T$ is weakly compact, then there exists a net $(T(e_j))_{j\in J}$ weakly converging to some $b_0\in B$. For every $a\in A$ we have
$$
T(e_j)\alpha(a)\rightarrow b_0\alpha(a)
$$
in the weak topology and 
$
T(e_j a)\rightarrow T(a)
$
in the norm topology. Hence $T(a)=\Lambda_{b_0, \alpha}(a)$. So (i) holds.

Assume now that $\alpha$ is weakly compact. Then there exists a net $(\alpha(e_j))_{j\in J}$ such that $\alpha(e_j)\rightarrow c_0\in B$ in the weak topology. Standard arguments show that $c_0$ is the identity of $\alpha(A)$. If $T$ is surjective, then $T(a_0)=c_0$ for some $a_0\in A$. So for every $a\in A$, we have
\begin{eqnarray*}
\alpha(a)&=&c_0\alpha(a)\\
&=&T(a_0)\alpha(a)\\
&=&T(a_0a)\\
&=&T\circ\lambda_{a_0}(a).
\end{eqnarray*}
Therefore (ii) holds.
$\hfill\square$

\section{Jordan left $\alpha$-centralizers on the group algebras} 

The existence of a bounded left approximate identity for $L^1(G)$ plays a pivotal role in analyzing linear operators. This structural feature leads to an important simplification: any continuous Jordan left $\alpha$-centralizer $T: L^1(G)\rightarrow L^1(G)$, where $\alpha\in\hbox{Hom}_c(L^1(G))$, is a left $\alpha$-centralizer.
This reduction from Jordan to ordinary centralizers provides the foundation for our next result in this section, which extends Wendel's celebrated characterization of multipliers [17] to the Jordan setting. 
%-----------------------------------------------------------------------
%-----------------------------------------------------------------------
%----------------------------------------------------------------------- 
\begin{theorem}\label{ame}
Let $G$ be a locally compact group and $\alpha\in\emph{Hom}_c(L^1(G))$.
If $T$ is a continuous Jordan left $\alpha$-centralizer on $L^1(G)$, then there exists a measure $\mu \in M(G)$ such that $T=\Lambda_{\mu, \alpha}$.
\end{theorem}
%-----------------------------------------------------------------------
%-----------------------------------------------------------------------
%-----------------------------------------------------------------------
{\it Proof.} 
Let $(e_i)$ be a bounded approximate identity for $L^1(G)$. Then $(T(e_i))_{i}$ forms a bounded net in $M(G)= C_0(G)^*$. 
Without loss of generality, we may assume that $T(e_i) \rightarrow \mu\in M(G)$ in the weak$^*$ topology of $M(G)$. 
So for every $\phi \in L^1(G)$, we have 
\begin{eqnarray*} 
T(e_i * \phi ) = T(e_i)* \alpha (\phi) \rightarrow \mu * \alpha (\phi)
\end{eqnarray*}
in the weak$^*$ topology of $M(G)$. 
On the other hand, $T(e_i *\phi) \rightarrow T(\phi)$ in the norm topology. Therefore, $T(\phi) = \mu * \alpha( \phi)=\Lambda_{\mu, \alpha}(\phi)$. 
$\hfill\square$\\

For $\alpha\in\hbox{Hom}(L^1(G))$, we denote by $\mathrm{JL}^c_\alpha(L^1(G))$ the set of all continuous Jordan left $\alpha$-centralizers on $L^1(G)$. It is easy to see that $\mathrm{JL}^c_\alpha(L^1(G))$ is a closed subspace of the Banach space of bounded linear operators on $L^1(G)$. So $\mathrm{JL}^c_\alpha(L^1(G))$ is itself a Banach space.

\begin{corollary}\label{mah} Let $\alpha\in\emph{Epm}_c(L^1(G))$. 
Then the mapping $\Gamma: M(G)\rightarrow\mathrm{JL}^c_\alpha(L^1(G))$ defined by $\Gamma(\mu)=\Lambda_{\mu, \alpha}
$ is a homeomorphism.
\end{corollary}
{\it Proof.} First, note that the map $\Gamma$ is surjective by Theorem 3.1. Let $(e_i)$ be a bounded approximate identity for $L^1(G)$ converging to $\delta_e$ in the weak$^*$ topology of $M(G)$, where $\delta_e$ is the Dirac measure at the identity element of $G$; see [6]. Since $\alpha$ is surjective, there exists a net $(\phi_i)_i$ in $L^1(G)$ such that $\alpha(\phi_i)=e_i$ for all $i$. To prove injectivity, let $\Gamma(\mu)=0$. Then 
$$
\mu\ast e_i=\mu\ast\alpha(\phi_i)=\Lambda_{\mu,\alpha}(\phi_i)=\Gamma(\mu)(\phi_i)=0.
$$
Taking weak$^*$ limits shows $\mu=0$. Thus $\Gamma$ is bijective and continuous, completing the proof by open mapping theorem. 
$\hfill\square$
%-----------------------------------------------------------------------
%-----------------------------------------------------------------------
%-----------------------------------------------------------------------
%-----------------------------------------------------------------------
%-----------------------------------------------------------------------
%-----------------------------------------------------------------------
\begin{theorem}\label{sm}
Let $G$ be a locally compact group, $\alpha\in\emph{Epm}_c(L^1(G))$ and $\mu\in M(G)$ be non-zero. Then the following assertions are equivalent.

\emph{(a)} $\Lambda_{\mu, \alpha}$ is compact on $L^1(G)$.

\emph{(b)} $\Lambda_{\mu, \alpha}$ is weakly compact on $L^1(G)$.

\emph{(c)} $G$ is compact and $\mu\in L^1(G)$.
\end{theorem}
%-----------------------------------------------------------------------
%-----------------------------------------------------------------------
%-----------------------------------------------------------------------
{\it Proof.}
Let $\mu\in M(G)$ be non-zero. Then $\lambda_\mu$ is non-zero on $L^1(G)$. Since $\alpha$ is surjective, it is routine to check that $\Lambda_{\mu, \alpha}$ is non-zero and (weakly) compact if and only if $\lambda_\mu$ is non-zero and (weakly) compact. By [2, 16], the left centralizer $\lambda_\mu$ is non-zero and compact if and only if $\lambda_\mu$ is non-zero and weakly compact which is equivalent to $G$ is compact and $\mu\in L^1(G)$. These observations complete the proof.
$\hfill\square$

%
%-----------------------------------------------------------------------
%-----------------------------------------------------------------------
%-----------------------------------------------------------------------
%-----------------------------------------------------------------------
%-----------------------------------------------------------------------
%-----------------------------------------------------------------------

\begin{remark}\label{hab0}
{\rm Let $G$ be a locally compact group and $\alpha\in\hbox{Epm}_c(L^1(G))$. Then the following statements hold.

(i) Every weakly compact Jordan left $\alpha$-centralizer $T$ on $L^1(G)$ admits a representation $T=\Lambda_{\phi, \alpha}$ for some $\phi\in L^1(G)$. 

(ii) The set of all weakly compact Jordan left $\alpha$-centralizers on $L^1(G)$ is either $\{0\}$ or homeomorphic to $L^1(G)$.

(iii) $L^1(G)$ admits a non-zero weakly compact Jordan left $\alpha$-centralizer if and only if for every $\phi\in L^1(G)$, the Jordan left $\alpha$-centralizer $\Lambda_{\phi, \alpha}$ is weakly compact; or equivalently, $G$ is compact.}
\end{remark}

In the following, we present some applications of Theorem 3.3.

\begin{theorem}\label{az}
Let $G$ be a locally compact group. Then the existence of 
a weakly compact epimorphism on $L^1(G)$ is equivalent to the finiteness of $G$.
\end{theorem}
{\it Proof.} Assume that $T$ is a weakly compact epimorphism on $L^1(G)$. Then $T$ is a weakly compact Jordan left $T$-centralizer on $L^1(G)$. According to Remark 3.4 (i), there exists $\psi\in L^1(G)$ such that $T=\Lambda_{\psi, T}$, that is, $T(\phi)=\psi\ast T(\phi)$ for all $\phi\in L^1(G)$. Let $(e_i)_i$ be a bounded approximate identity for $L^1(G)$ such that $e_i\rightarrow\delta_e$ in the weak$^*$ topology of $M(G)$. Since $T$ is surjective, there exists a net $(\phi)_i$ in $L^1(G)$ with $T(\phi_i)= e_i$ for all $i$. Hence
$$
e_i=T(\phi_i)=\psi\ast T(\phi_i)=\psi\ast e_i\rightarrow\psi
$$
in the norm topology. This implies that $\delta_e=\psi\in L^1(G)$. Hence $G$ must be discrete. Applying Remark 3.4 (iii), we conclude that $G$ is finite.
The converse is straightforward.
$\hfill\square$

\begin{corollary}\label{do}
Let $G$ be a locally compact group and $\alpha\in\emph{Epm}_c(L^1(G))$. Then every continuous Jordan left $\alpha$-centralizer on $L^1(G)$ is weakly compact if and only if $G$ is finite.
\end{corollary}

Let $\alpha\in\hbox{Hom}(L^1(G))$. Then the linear map $D: L^1(G)\rightarrow L^1(G)$ is called an $\alpha$-\emph{derivation} if
$$
D(\phi\ast\psi)=D(\phi)\ast\alpha(\psi)+\alpha(\phi)\ast D(\psi)
$$
for all $\phi, \psi\in L^1(G)$.

\begin{theorem}\label{dev} Let $G$ be a locally compact group and $\alpha\in\emph{Epm}_c(L^1(G))$. Then $L^1(G)$ admits a non-zero weakly compact $\alpha$-derivation if and only if $G$ is compact and non-abelian.
\end{theorem}
{\it Proof.} Let $D: L^1(G)\rightarrow L^1(G)$ be a non-zero weakly compact $\alpha$-derivation. For every $\phi, \psi\in L^1(G)$, we have
\begin{eqnarray*}
\Lambda_{D(\phi), \alpha}(\psi)&=&D(\phi)\ast\alpha(\psi)\\
&=&D(\phi\ast\psi)-\alpha(\phi)\ast D(\psi)\\
&=&D\lambda_\phi(\psi)-\lambda_{\alpha(\phi)}D(\psi).
\end{eqnarray*}
This shows that $\Lambda_{D(\phi), \alpha}$ is weakly compact. Choose $\phi_0\in L^1(G)$ such that $D(\phi_0)$ is non-zero. Since $\alpha$ is surjective, $\Lambda_{D(\phi_0), \alpha}$ is non-zero. Hence $L^1(G)$ admits a non-zero weakly compact Jordan left $\alpha$-centralizer. Therefore, $G$ is compact.

 Suppose now that $G$ is abelian. Then $L^1(G)$  is a commutative semisimple Banach algebra. So the Gel'fand transform $$\Gamma: L^1(G)\rightarrow C_0(\Delta(L^1(G)))$$ is a monomorphism [3, 4], where $\Delta(L^1(G))$ is the spectrum of $L^1(G)$. For $\chi\in\Delta(L^1(G))$, define $d_\chi: L^1(G)\rightarrow{\Bbb C}$ and $\tilde{\chi}: L^1(G)\rightarrow{\Bbb C}$ by 
$$
d_\chi(\phi)=\chi(D(\phi)\quad\hbox{and}\quad\tilde{\chi}(\phi)=\chi(\alpha(\phi)).
$$
It is easy to see that $\tilde{\chi}\in \Delta(L^1(G))$ and $d_\chi$ is a point derivation at $\tilde{\chi}$; that is,
$$
d_\chi(\phi\ast\psi)=d_\chi(\phi)\ast\tilde{\chi}(\psi)+\tilde{\chi}(\phi)\ast d_\chi(\psi)
$$
for all $\phi, \psi\in L^1(G)$. Since $L^1(G)$ is weakly amenable [15] and no weakly amenable Banach algebra admits non-zero continuous point derivations [4], it follows that $\chi(D(\phi))=0$ for all $\phi\in L^1(G)$. As $\chi$ is arbitrary, $\Gamma(D(\phi))=0$ for all $\phi\in L^1(G)$. Thus $D(\phi)=0$ for all $\phi\in L^1(G)$, contradicting  the non-zero assumption. Therefore, $G$ must be non-abelian.

For the converse, let $G$ be compact and non-abelian. Then $L^1(G)$ is non-commutative. So there exist $\phi, \psi\in L^1(G)$ with $$\phi\ast\psi\neq\psi\ast\phi$$ and $\phi\neq 0$. It follows from Theorem 3.3 that $\Lambda_{\phi, \alpha}$ is weakly compact. Likewise, an argument analogous to the proof of Theorem 3.3 shows that the mapping $\Sigma_{\alpha, \phi}: L^1(G)\rightarrow L^1(G)$ defined by $$\Sigma_{\alpha, \phi}(\eta)=\alpha(\eta)\ast\phi$$ is weakly compact. Define $D: L^1(G)\rightarrow L^1(G)$ as $$D=\Lambda_{\phi, \alpha}-\Sigma_{\alpha, \phi}.$$
Then $D$ is a weakly compact $\alpha$-derivation
on $L^1(G)$. By surjectivity of $\alpha$, there exists $\psi_0\in L^1(G)$ such that
$\alpha(\psi_0)=\psi$. Then 
\begin{eqnarray*}
D(\psi_0)&=&\Lambda_{\phi, \alpha}(\psi_0)- \Sigma_{\alpha, \phi}(\psi_0)\\
&=&\phi\ast\alpha(\psi_0)-\alpha(\psi_0)\ast\phi\\
&\neq& 0,
\end{eqnarray*}
proving $D$ is non-zero.
$\hfill\square$\\

In Corollary 3.2, Theorems 3.3 and 3.5, Corollary 3.6, and Theorem 3.7 we assume that the homomorphism $\alpha$ is surjective. It is natural to ask whether These results hold without this assumption. We present counterexamples demonstrating that this condition can not be remove. First, we note that if $G$ is compact and $\mu\in L^1(G)$, then the left centralizer $\lambda_\mu$ is (weakly) compact on $L^1(G)$; see [2]. Therefore, $\Lambda_{\mu, \alpha}$ is also (weakly) compact on $L^1(G)$. This proves that the ``if" direction of Theorem 3.3 remains valid even without assuming $\alpha$ is surjective. However, as the following examples show, the surjectivity of $\alpha$ is essential for the ``only if" direction of the theorem.

\begin{example}\label{224}{\rm (i) Let $G={\Bbb R}$, the additive group of real numbers with the usual topology. Define the homomorphism $\alpha: L^1(G)\rightarrow L^1(G)$ by $\alpha(\phi)=0$ for all $\phi\in L^1(G)$. Note that $\alpha$ is not surjective and for every $\mu\in M(G)$, we have $\Lambda_{\mu, \alpha}=0$ on $L^1(G)$. Thus
$$
\mathrm{JL}^c_\alpha(L^1(G))=\{0\}
$$
and so every continuous Jordan left $\alpha$-centralizer on $L^1(G)$ is weakly compact. Although $G$ is a non-compact infinite group, $\Lambda_{\delta_0, \alpha}$ is trivially weakly compact on $L^1(G)$, yet $\delta_0\not\in L^1(G)$. Hence Corollary 3.2, the ``only if" direction of Theorem 3.3, Theorem 3.5 and Corollary 3.6 fail without the surjectivity assumption on $\alpha$.

(ii) Let $G$ be the circle group ${\Bbb T}={\Bbb R}/{\Bbb Z}$. Define $\alpha: L^1({\Bbb T})\rightarrow L^1({\Bbb T})$ by 
$$
\alpha(\phi)(x)=\int_{\Bbb T}\phi d\lambda\quad\quad(\phi\in L^1(G), x\in {\Bbb T}),
$$
where $\lambda$ is the normalized Haar measure on ${\Bbb T}$. Thus $\alpha(\phi)$ is the constant function equal to $\int_{\Bbb T} \phi d\lambda$. Then $\alpha$ is a homomorphism but it is not surjective. 

For every $\mu\in M(G)$ and $\phi\in L^1(G)$, we have
$$
\mu\ast\alpha(\phi)=(\int_{\Bbb T}\phi d\lambda)\mu(G).
$$
Hence $\Lambda_{\delta_s, \alpha}=\Lambda_{\delta_t, \alpha}$ for all $s, t\in{\Bbb T}$. This implies that $\Gamma(\delta_s)=\Gamma(\delta_t)$ for distinct elements $s, t\in{\Bbb T}$. This shows that $\Gamma$ is not injective.

Since $\alpha$ is a finite-rank operator, it is weakly compact. Hence every continuous Jordan left $\alpha$-centralizer on $L^1(G)$ is weakly compact. In particular, $\Lambda_{\delta_1, \alpha}$ is weakly compact, though $\delta_1\not\in L^1({\Bbb T})$ and ${\Bbb T}$ is infinite.

This example confirms the necessity of surjectivity in Corollary 3.2, the ``only if" direction of Theorem 3.3, Theorem 3.5 and Corollary 3.6.

(iii) Let \( G = \hbox{SO}(3) \), the compact non-abelian group of $3D$ rotations. Let $\alpha: L^1(G) \to L^1(G)$ be the constant function with the value
$$
\alpha(\phi)(x) = \int_G \phi  d\lambda,
$$
where $\lambda$ is the normalized Haar measure on $G$. Then $\alpha$ is homomorphism but not surjective. Define the multiplicative linear functional $\chi$ on $L^1(G)$ by $\chi(\phi)= \int_G \phi  d\lambda$. Note that for every $\phi, \psi\in L^1(G)$ we obtain
\[
\phi \ast\alpha(\psi) =\chi(\psi)\int_G \phi  d\lambda= \alpha(\psi)\ast \phi.
\] 
Let \( D \) be a weakly compact \( \alpha \)-derivation. Then for every $\phi, \psi\in L^1(G)$ we have
\begin{eqnarray}\label{er}
D(\phi \ast \psi)&=& D(\phi) \ast \alpha(\psi) + \alpha(\phi) \ast D(\psi)\nonumber\\
&=&
\chi(\psi) \int_G D(\phi)  d\lambda+ \chi(\phi)\int_G D(\psi)  d\lambda.
\end{eqnarray}
Define the linear map $d: L^1(G)\rightarrow{\Bbb C}$ by
$
d(\phi)=\int_G D(\phi) d \lambda.
$
From (6) we infer that
$$
d(\phi\ast\psi)= d(\phi)\chi(\psi)+\chi(\phi) d(\psi)
$$
for all $\phi, \psi\in L^1(G)$. That is, $d$ is a point derivation at $\chi$. Since $G$ is compact, $L^1(G)$ is amenable and so it is weakly amenable. Thus, $d=0$. Therefore $D=0$. Consequently, the only $\alpha$-derivation on $L^1(G)$ is zero. This shows that surjectivity of \( \alpha \) is necessary for Theorem 3.7.
}
\end{example}

\vspace{2mm}

{\footnotesize
\noindent {\bf Mojdeh Eisaei}\\
Department of Mathematics, \\Payame Noor University,\\ Shiraz, Iran\\ e-mail: mojdehessaei59@student.pnu.ac.ir\\
\noindent {\bf Mohammad Javad Mehdipour}\\
Department of Mathematics,\\ Shiraz University of Technology,\\
Shiraz
71555-313, Iran\\ e-mail: mehdipour@sutech.ac.ir\\
\noindent {\bf Gholam Reza Moghimi}\\
Department of Mathematics, \\Payame Noor University,\\ Shiraz, Iran\\ e-mail: moghimimath@pnu.ac.ir\\

\footnotesize

\begin{thebibliography}{99}

\bibitem{am} M. H. Ahmadi Gandomani and M. J. Mehdipour, Symmetric bi-derivations and their generalizations on group algebras. Filomat 35 (2021), no. 4, 1233--1240.

\bibitem{a} C. A. Akemann, Some mapping properties of the group algebras of a compact group. Pacific J. Math. 22 (1967), 1--8.

\bibitem{bd} F. F. Bonsall and J. Duncan, Complete Normed Algebras, Springer Verlage, Berlin, Heidelberg and New York, 1973.

\bibitem{d} H. G. Dales, Banach algebras and automatic continuity, London Math. Soc. Monographs 24, Oxford Univ. Press, New
York, 2000.

\bibitem{g} F. Ghahramani, Compact elements of weighted group algebras. Pacific J. Math. 113 (1984), no. 1, 77--84. 

\bibitem{g1} F. Ghahramani, Weighted group algebra as an ideal in its second dual space. Proc. Amer. Math. Soc. 90 (1984), no. 1, 71--76.

\bibitem{gl} F. Ghahramani and A. T. Lau, Isomorphisms and multipliers on second dual algebras of Banach algebras. Math. Proc. Cambridge Philos. Soc. 111 (1992), no. 1, 161--168.

\bibitem{gm} M. Ghasemi and M. J. Mehdipour, Homological properties of Banach modules related to locally compact groups. Proc. Rom. Acad. Ser. A Math. Phys. Tech. Sci. Inf. Sci. 21 (2020), no. 4, 295--301

\bibitem{hr} E. Hewitt and K. Ross, Abstract Harmonic Analysis I, Springer, New York, 1970.

\bibitem{lp} A. T. Lau and J. Pym, Concerning the second dual of the group algebra of a locally compact group. J. London Math. Soc. (2) 41 (1990), no. 3, 445--460.

\bibitem{l} V. Losert, Weakly compact multipliers on group algebras. J. Funct. Anal. 213 (2004), no. 2, 466--472.

\bibitem{m} D. Malekzadeh Varnosfaderani, Derivations, multipliers and topological centers of certain Banach algebras related to
locally compact groups, Ph.D. thesis, University of Manitoba, 2017, available at https://mspace.lib.umanitoba.ca/xmlui/
handle/1993/32276.

\bibitem{mm} M. J. Mehdipour and Gh. R. Moghimi, The existence of nonzero compact right multipliers and Arens regularity of weighted Banach algebras. Rocky Mountain J. Math. 52 (2022), no. 6, 2101--2112. 

\bibitem{mn2} M. J. Mehdipour and A. Rejali, Regularity and amenability of weighted Banach algebras and their second dual on locally compact groups, arXiv:2112.13286

\bibitem{r} V. Runde, Lectures on Amenability, Springer Verlage, Berlin, Heidelberg, 2002.

\bibitem{s} S. Sakai, Weakly compact operators on operator algebras. Pacific J. Math. 14 (1964), 659--664.

\bibitem{w} J. G. Wendel, Left centralizers and isomorphisms of group algebras. Pacific J. Math. 2 (1952), 251--261

\bibitem{z} B. Zalar, On centralizers of semiprime rings. Comment. Math. Univ. Carolin. 32 (1991), no. 4, 609--614.

\end{thebibliography}
\end{document}